# Multi-Objective Optimization of a Port-of-Entry Inspection Policy

Christina M. Young, Mingyu Li, Yada Zhu, Minge Xie, Elsayed A. Elsayed, and Tsvetan Asamov

E. A. Elsayed, C. M. Young, and Y. Zhu are with the Department of Industrial and Systems Engineering, Rutgers University, Piscataway, NJ 08854 USA (e-mail: elsayed@rci.rutgers.edu). M. Xie and M. Li are with the Department of Statistics, Rutgers University, Piscataway NJ 08854 USA (e-mail: mxie@stat.rutgers.edu). T. Asamov is with Kenyon College

**Abstract**

At the port-of-entry containers are inspected through a specific sequence of sensor stations to detect the presence of nuclear materials, biological and chemical agents, and other illegal cargo. The inspection policy, which includes the sequence in which sensors are applied and the threshold levels used at the inspection stations, affects the probability of misclassifying a container as well as the cost and time spent in inspection. In this paper we consider a system operating with a Boolean decision function combining station results and present a multi-objective optimization approach to determine the optimal sensor arrangement and threshold levels while considering cost and time. The total cost includes cost incurred by misclassification errors and the total expected cost of inspection, while the time represents the total expected time a container spends in the inspection system. An example which applies the approach in a theoretical inspection system is presented.

## I. INTRODUCTION

The significant increase in trade agreements and the growth in the world economy have propelled unprecedented increase in maritime traffic. The value of export goods produced and transported globally in 2000 was about $6.186 trillion [1]. Disruption of such a system has catastrophic consequences on the world economy and our daily needs. In order to minimize sources of disruptions, the United Nations passed several resolutions with the objective of improving security in maritime trade. Likewise, the United States initiated the Container Security Initiative (CSI) to ensure container security through different approaches starting from the origin port of the container and ending at the delivery port in the United States. When containers arrive at United States ports they can be randomly selected and subjected to inspection. The type of inspections, number of containers to be inspected, and the inspection policy have a profound effect on the cost of the system, risk of accepting undesired containers and potential delays and congestion at the ports.

In this paper we consider a port-of-entry (POE) container inspection system where a fraction of the arriving containers at a port is subjected to a sequence of inspections at different stations. A typical inspection system involves radiation detection and may include x-ray or gamma-ray imaging, as well as methods currently under research such as biosensors and chemical agent detection.

Researchers have investigated the problem of container inspection with different objectives. Lewis *et al.* [2] develop a best-first heuristic search procedure to model the problem of moving containers to and from inspection areas, but the inspection procedures and sequences have not been considered. Stroud and Saeger [3] consider the problem where containers arrive at a port and sequential inspections are conducted to decide whether to pass a container or subject it to further inspection. Containers that meet some conditions and leave the inspection system can either be accepted or subjected to "manual" inspection. This problem is considered as a binary decision tree problem. Madigan *et al.* [4] extend the work of Stroud and Saeger by incorporating the threshold levels of the inspection "sensors" and develop a novel binary decision tree search algorithm that operates on a space of potentially acceptable binary decision trees. They describe computationally more efficient approaches for this binary decision tree problem and obtain optimum sensor threshold levels that minimize the total cost of the inspection system.

Elsayed *et al.* [5] present a unique approach to the formulation of the port-of-entry inspection problem as an analytical model. Unlike previous work which determines threshold levels and sequence separately, they consider an integrated system and determine them simultaneously. They decompose the POE problem into two sub-problems. One problem deals with the determination of the optimum sequence of inspection or the structure of the inspection decision tree in order to achieve the minimum expected inspection cost. This problem is formulated and investigated using approaches parallel to those used in the optimal sequential inspection procedure for reliability systems as described in [6]-[13]. The other problem deals with the determination of the optimum thresholds of the inspection stations so as to minimize the cost associated with false reject and false accept. As indicated earlier, the delay in inspection system is also a major concern as it has significant economic consequences.

In this paper we develop a new formulation for the POE problem by considering both the minimization of the total cost and the delay time of the containers simultaneously as a multi-objective optimization problem. We seek the optimum inspection sequence and the optimum threshold levels of sensors at inspection stations in order to minimize the total cost and total delay time.

## II. PROBLEM DESCRIPTION

### A. Port-of-Entry Container Inspection System

In modeling the port-of-entry container inspection system it is assumed that containers arriving for inspection are inherently acceptable or contain unacceptable materials, and that they have several attributes which may reflect the status (presence or no presence) of such material. The inspection system is viewed as a collection of stations, over which the inspection of a given container is performed sequentially. Each station inspects one specific attribute and returns a pass-or-fail decision (0 or 1 respectively). At each individual station the decision is dependent on a preset threshold level. Varying this threshold level affects the probability of misclassifying an acceptable container as suspicious or vice versa. The sequence in which stations are to be visited, along with threshold levels to be applied, establishes the inspection policy which is applied to every container arriving for inspection. The final decision to accept a container or reject it, thereby subjecting the container to further manual inspection, often including a manual "unpacking" method, is determined based on the evaluation of a predefined Boolean decision function of the individual station decisions.

The Boolean decision function $F$ assigns to each binary string of attributes $(a_1, a_2, ..., a_n)$ a category. In other words $F(a_1, a_2, ..., a_n) = 0$ indicates negative class and that there is no suspicion with the container and $F(a_1, a_2, ..., a_n) = 1$ indicates positive class and that additional inspection is required, usually manual inspection.

By definition, for instance, a series Boolean function is a decision function $F$ that assigns the container a value of "1" if any of the attributes is present, i.e. $a_i = 1$ for any $i \in \{1, 2, ..., n\}$, and a parallel Boolean function is a decision function $F$ that assigns the container a value of "1" if all of the attributes are present, i.e. $a_i = 1$ for all $i \in \{1, 2, ..., n\}$. The Boolean function to be used depends on the nature of the inspection system; the container attributes being inspected, and other factors. The work presented here is designed so that it can be applied with any Boolean function.

### B. Modeling of Sensor Measurements and Thresholds

Let $x$ represent true status of a container, and code $x = 1$ if it should be rejected and $x = 0$ if it should be accepted. We assume this container is a sample from a population of interest under which the probability of $x = 1$ is $P(x = 1) = \pi$ and the probability of $x = 0$ is $P(x = 0) = 1 - \pi$.

Let $r$ be the measurement taken by a sensor. This measurement $r$ can in general be a numerical (continuous or discrete) reading or a graphical image. To simplify the presentation of our development and following [3] and [5], we assume $r \sim N(\mu_0, \sigma_0^2)$ when $x = 0$ and $r \sim N(\mu_1, \sigma_1^2)$ when $x = 1$, where $\mu_0 \neq \mu_1$. We choose to use the normal distribution because normally distributed data are the most commonly seen data in practice and it is has been used in port-of-entry inspection applications [3], [5]. Also, continuous measurements can often be transformed into a normal distributed data by the well known inverse transformation method [14]. Likewise, discrete data sometimes can be well approximated by a normal distribution either by central limit theorem or by some special techniques such as variance stabilization transformation. Our development, in principle, can be extended to some non-normal cases.

We assume two normal distributions for $r$ because we expect to have different sensor readings for a container with true status $x = 1$ and $x = 0$. We also assume that the parameters of the two normal distributions indicated are known or can be estimated from past inspection history. Note that the task of distinguishing acceptable and unacceptable containers is location and scale invariant to the readings. Without loss of generality, we can assume that $\mu_0 = 0$ and $\mu_1 = 1$. See also [5] for further discussions on this assumption.

To make a decision based on the sensor measurement $r$, the $r$ value is compared against a given threshold value $T$. We reject the container ($d = 1$) if the reading $r$ is higher than $T$ and accept it ($d = 0$) if the reading is less than or equal to $T$. The decision $d$ at this level of decision making does not always agree with the true status $x$. There are two types of potential errors. A type I error is a decision $d = 1$ when the true status of the container is $x = 0$, and a type II error is a decision $d = 0$ when the true status of the container is $x = 1$. The probability of these two types of errors can be computed by 
$$P(d = 1 | x = 0) = P(r > T | x = 0) = 1 - \Phi\left(\frac{T}{\sigma_0}\right)$$ 
and
$$P(d = 0 | x = 1) = P(r \leq T | x = 1) = \Phi\left(\frac{T - 1}{\sigma_1}\right).$$

### C. System Inspection Policy

The minimization of costs associated with performing inspection and misclassification of containers has been formulated in Elsayed *et al.* [5]. Here we expand the optimization objective to include the time required for inspection, which takes into account the effect of delays on the overall system. The time incurred in inspection is added to the objectives because it may be very important in some situations. The performance of the inspection system is determined by both the sequence in which inspection stations are visited and the threshold levels applied at those stations, which we denote collectively as the inspection policy.

Since the optimal parameter values for the cost minimization problem may not minimize time, some compromise may be required. A particular balance of the importance of cost and time may be represented by weights. We consider the case where the relative importance of cost and time is unspecified and therefore we use different importance weights to generate possible solutions that produce a Pareto frontier as described in section IV.

### III. PERFORMANCE MEASURES OF INSPECTION POLICY

### A. Cost of Misclassification and Inspection

The cost involved in this inspection problem is the sum of any cost incurred as a result of misclassifying a container's status and the actual cost of performing the inspection. As Elsayed *et al.* [5] note, there are two types of misclassification errors at the systems level: falsely rejecting a container that should be cleared and falsely accepting a container that should be rejected. These errors are associated with the probability of false reject (PFR) and the probability of false accept (PFA), respectively. The complementary probabilities of these two errors are true reject (PTR) and true accept (PTA). If $D$ denotes the decision of the entire inspection system of sensors, where $D = 1$ means to reject and $D = 0$ to accept, the four probabilities can be written as follows:

$PFR = P(D = 1 | x = 0)$, $PTA = P(D = 0 | x = 0) = 1 - PFR$,

$PFA = P(D = 0 | x = 1)$, and $PTR = P(D = 1 | x = 1) = 1 - PFA$.

The inspection decision $D$ depends on the individual inspection results and the system Boolean function. The probability equations just mentioned can be rewritten in terms of the threshold value $T_i$ and $\sigma$ values related to the inspection station for any given Boolean function. Several examples are given in Elsayed *et al.* [5].

The cost of misclassification arises when a cost is associated with PFR and PFA. Let $c_{FA}$ be the cost of the system accepting a "bad" container and $c_{FR}$ be the cost of the system rejecting a "good" container. Then the total cost of system misclassification error is $C_F = \pi\, PFA\, c_{FA} + (1-\pi) PFR\, c_{FR}$ as described in [5].

The expectation of the cost of inspection is a function of the unit cost to operate each sensor (station) and the probabilities of passing each station. Given a particular set of threshold values, an optimal sequence in which to visit the sequence can be found following the conditions in Elsayed *et al.* [5]. The total cost arising from misclassification errors and inspection is denoted by $c_{total} = C_F + E[C_{inspection}]$.

### B. Time for Inspection

The time required for a container to pass through an inspection station is an important measure of the inspection system performance. It is possible that this time would be related to some characteristic of the inspection station, that is to say the inspection may be sped up or slowed down depending on some operational setting of the sensor. For example the inspection time may be related to a variable that represents the resolution or other settings of the sensor. Following Jupp *et al.* [15], we assume the time spent at each station could be related to the threshold level $T_i$ at that station. This relationship is expressed as $t_i = a\exp(b \cdot T_i)$ for illustration purposes. Here the time of inspection decreases as the applied threshold level increases.

To find the total expectation of time spent in the system for a given container we first denote $p_i$, the probability of passing station $i$, by:

$$p_i = P(d_i = 0) = \sum_{j=0}^{1}\left[P(d_i=0\,|\,x=j)P(x=j)\right],$$ which can

be written $p_i = (1-\pi)\Phi\left(\dfrac{T_i}{\sigma_{0i}}\right) + \pi\Phi\left(\dfrac{T_i-1}{\sigma_{1i}}\right)$, and $q_i = 1 - p_i$,

where $p_i$ and $q_i$ are functions of threshold values $T_i$. Then, the total expected inspection time $t_{total}$ can be expressed as $t_{total} = t_1 + \sum_{i=2}^{n}\left[\prod_{j=1}^{i-1} p_j\right] t_i$, where $t_i$ is the inspection time at station $i$, for a system with $n$ stations using a series Boolean decision function. For a parallel Boolean decision function, the total expected inspection time is $t_{total} = t_1 + \sum_{i=2}^{n}\left[\prod_{j=1}^{i-1} q_j\right] t_i$.

## IV. MULTI-OBJECTIVE OPTIMIZATION

### A. Total Expected Cost and Time

As noted in the problem description, we need to determine the optimal configuration of sensors in the system and the optimum sets of threshold levels that can achieve the objectives of maximizing inspection system throughput and minimizing the expected total cost. The POE problem formulation can be expressed as: $\min_{Sequence, Threshold}\{c_{total}, t_{total}\}$. It is unlikely that these objectives would be optimized by the same set of parameter values, and there exists some trade-off between them. This is a typical multi-objective optimization problem. See, for instance, Eschenauer *et al.* [16], Statnikov and Matuso [17], Fonseca and Fleming [18], [19], and Leung and Wang [20], among others.

The multi-objective problem is often solved by combining the multiple objectives into one scalar objective whose solution is a Pareto optimal point for the original problem. In general, there may be a large number or infinite number of optimal solutions in the sense of Pareto-optimality. In the POE case it is desirable to find as many optimal solutions as possible in order to provide more choices to decision makers. A commonly used method is the weighted sum approach, where weighted sums of the objective functions are optimized for various choices of fixed weights $w_1$ and $w_2$, $w_1 + w_2 = 1$. The fitness function used in this work is $f_{w_1,w_2}(S,T) = w_1 c_{total} + w_2 t_{total}$. Here, $S$ and $T$ stand for sequence and threshold levels that comprise an inspection policy. Thus, the multi-objective optimization problem becomes a collection of single objective optimization problems, in which we minimize the fitness function for a set of fixed weights $w_1$ and $w_2$: $\min_{S,T} f_{w_1,w_2}(S,T)$.

It is computationally expensive to directly solve this minimization problem because the number of potential sequences grows exponentially as the number of inspection stations increases. Therefore we employ a modified weighted sum approach, in which we utilize some theoretical results to decide the system sequence. Note that the fitness function is highly discrete with regards to the system sequence.

For the system Boolean functions considered in this paper, the optimal sequence can be obtained for a given set of weights and thresholds. Theorem 1 in [5] demonstrates that an optimal sequence for a series Boolean system can be formed by arranging the sensors $i = 1,...,n$ by the value of the ratio $c_i/q_i$, from least to greatest. We extend this to the case of weighted multiple objectives by using $(w_1 c_i + w_2 t_i)/q_i$ as the ratio.

The value of the fitness function can be calculated by:

$$f_{w_1,w_2}(T) = (w_1 c_1 + w_2 t_1) + \sum_{i=2}^{n}\left[\prod_{j=1}^{i-1} p_j\right](w_1 c_i + w_2 t_i) + w_1 C_F \quad \text{and}$$

this expression is minimized when the stations are arranged in the optimal sequence for a series Boolean system. Likewise, the optimal inspection sequence can be found for a parallel Boolean using the ratio $c_i/p_i$, and the value of the fitness function can be written by substituting $q_i$ in place of $p_i$. The optimal sequence theory for series system and parallel system can be extended to systems with arrangements of parallel-series and series-parallel sensors similarly to how Theorem 2 is stated in [5], by simply replacing all instances of $c_i$ with $w_1 c_i + w_2 t_i$. So, for a given set of thresholds and specified weights the optimal inspection sequence can be found and allows us to compute the function $f_{w_1,w_2}(T) = \min_{S} f_{w_1,w_2}(S,T)$.

It then remains to solve the minimization problem $\min_{T} f_{w_1,w_2}(T)$. This method can be easily modified to apply to systems using many different Boolean decision functions. This modified approach can provide an efficient method to deal with the multi-objective optimization problem under the current context by avoiding the consideration of all possible sequences.

*B. Implementation: Three Approaches*

Three methods are distinguished in the implementation: Grid Search (GS), *fmincon* (FM) and Genetic Algorithm (GA). The grid search method is a complete enumeration method. It sets a standard against which the GA and FM approaches may be compared. FM and GA are based on the optimization algorithm developed in the previous subsection. The difference between these methods is how they solve the optimization problem $\min_T f_{w_1,w_2}(T)$.

The grid search method is a complete enumeration of possible threshold values and all inspection sequences. A discrete set of threshold values is formed from the range of 0-1 using a gradient of 0.05. The total cost and total time are calculated for each possible combination of threshold values and sequence. The resulting cost and time values are plotted and the outermost points along the curve are filtered to represent the solution set that forms the Pareto frontier. Thus the GS method yields a small number of true optimal points compared to the other methods.

The MATLAB (The MathWorks, Inc.) *fmincon* function attempts to find a constrained minimum of scalar function of several variables starting at an initial estimate. This is generally referred to as constrained nonlinear optimization or nonlinear programming. For each pair of weights, we use *fmincon* to minimize $f_{w_1,w_2}(T)$ and try different sets of initial thresholds. Note that the function $f_{w_1,w_2}(T) = \min_S f_{w_1,w_2}(S,T)$ is highly discrete in $T$ which is inherited from the sequence optimization. Therefore, direct use of the *fmincon* function does not always work. As the third subplot of Figure 1 shows, the optimal solutions vary significantly with different initial values used in the *fmincon* function.

Finally, a genetic algorithm is an iterative random search algorithm, which takes advantage of information in the previous steps (ancestors) to produce new searching points (off-springs). It is called "genetic" algorithm because the principle and design of this search algorithm mimics those of genetic evolution found in nature [21]. A genetic algorithm can be applied to solve "a variety of optimization problems that are not well suited for standard optimization algorithms, including problems in which the objective function is discontinuous, nondifferentiable, stochastic, or highly nonlinear" [21]. In the optimization algorithm developed here, the MATLAB function *ga* was used to minimize $f_{w_1,w_2}(T)$ for each pair of weights. One advantage of using the *ga* is that it is insensitive to the initial values and we are able to obtain optimal solutions in all of our analyses.

## V. NUMERICAL EXAMPLE

Here the results of the multi-objective optimization are presented in graphical form. The three graphs in Figure 1 illustrate the optimal points obtained from the three methods discussed in the previous section applied to an inspection system using a parallel Boolean decision function. The system parameters in this example are as follows: $n=3$, $c=[1\ 1\ 1]$, $\pi=0.0002$, $\mu_0=[0\ 0\ 0]$, $\mu_1=[1\ 1\ 1]$, $\sigma_0=[0.16\ 0.2\ 0.22]$, $\sigma_1=[0.3\ 0.2\ 0.26]$, $c_{FA}=100000$, $c_{FR}=500$, $a=[20\ 20\ 20]$, $b=[-3\ -3\ -3]$, $w_1=[0:\ 0.004:1]$, $w_2=1-w_1$.

The grid search method produces optimal points that fall into distinct vertical segments due to the discrete nature of the method, and the minimum search gradient with an acceptable computation time was used. The leftmost graph contains only the outermost points with respect to the Pareto frontier from this method. Note that a small number of the points shown actually fall on the theoretical Pareto frontier, therefore the output from this method is not as useful compared to the others.

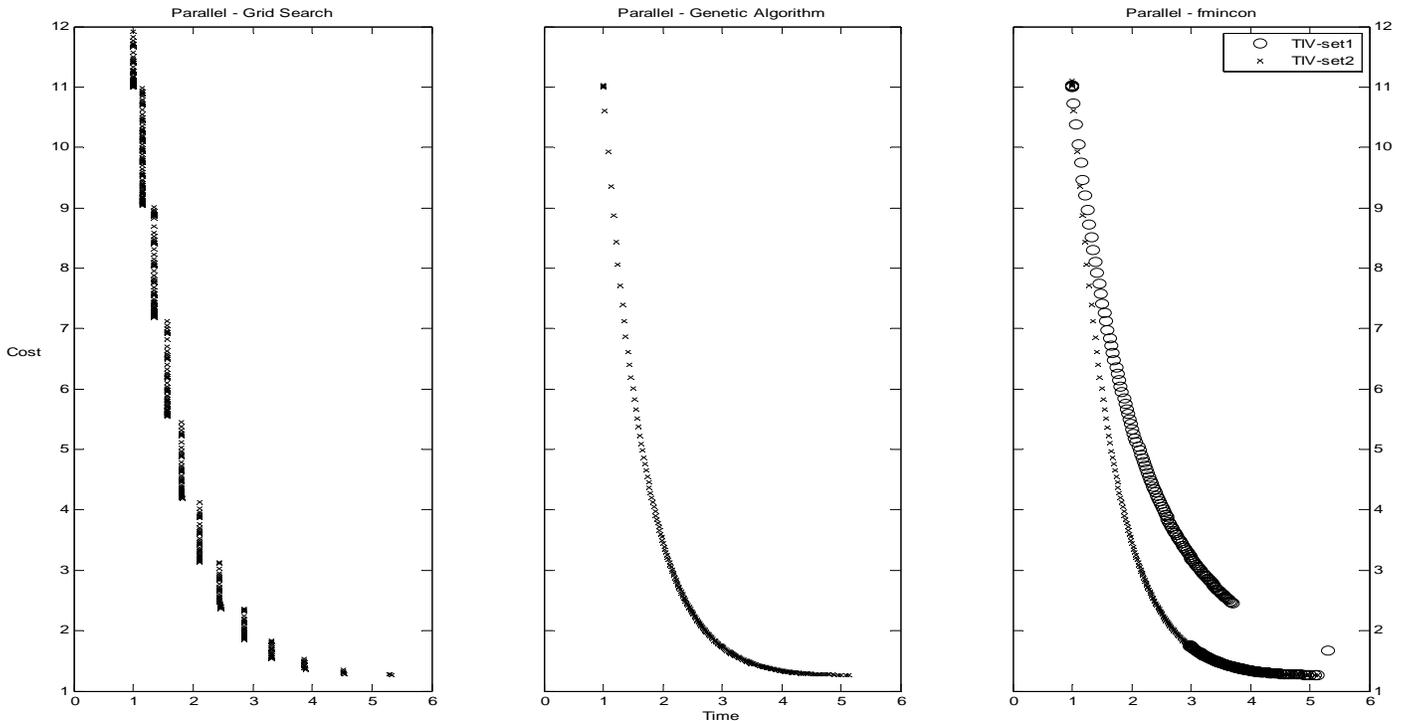

**Figure 1. Comparison of Three Solution Methods to Multi-Objective Problem**

The center graph illustrates the optimal points obtained from the GA method. For the FM method it was discovered that the initial values used had a significant effect on the optimality of the points obtained. Therefore various sets of initial values were tested and the results overlaid on one graph to illustrate the phenomenon. The initial value sets are represented as two series in the rightmost graph. The initial values for threshold initial value (TIV) set 1 are $T_0$=[0.2 0.2 0.2] and the initial values for TIV-set2 are $T_0$=[0.2 0.6 0.2]. Another set, $T_0$=[0.8 0.8 0.8] gave similar results to TIV-set1, and is not shown.

All three methods produce at least some portion of the same Pareto frontier of solutions with minimal time and cost. Each point represents the time and cost for one possible solution, and each solution is defined by a set of threshold values $\{T_i : i = 1,...,n\}$ -each to be applied at one of the $n$ inspection stations- and the sequence in which to visit those stations. Table 1 presents three examples of points chosen from the Pareto frontier of the grid search results.

**Table 1. Examples of Pareto Optimal Solutions**

| T1 | T2 | T3 | Sequence | Cost | Time |
|---|---|---|---|---|---|
| 0.0 | 0.95 | 0.05 | 2-3-1 | 9.03 | 1.16 |
| 0.0 | 0.85 | 0.0 | 2-1-3 | 5.54 | 1.57 |
| 0.0 | 0.75 | 0.05 | 2-3-1 | 3.13 | 2.11 |

It is important to consider program running time in the comparison of methods. FM is the fastest of the three methods, requiring about 2 minutes for one initial set. However, it does not always return the correct Pareto frontier and thus different initial sets must be used, and without knowledge of the true Pareto frontier choosing a good initial set is difficult. The GS method with grid=0.05 runs in about 6 minutes, however only about 12 points of the output are considered to fall within the theoretical Pareto frontier. If the grid is decreased to 0.025, roughly 23 points on the theoretical Pareto frontier are produced but it takes 5 hours to run. Further reducing the grid to 0.01 requires more than 200 hours to finish. Therefore it becomes impractical to decrease the grid size in order to generate more optimal points on the theoretical frontier.

The GA method takes about 10.5 hours with the current choice of parameters (PopulationSize=80) and produces 251 points on the theoretical Pareto frontier. Note that the *ga* function of MATLAB is designed for general purpose use, and we anticipate that the running time can be significantly improved by using a specialized program. Moreover, the GA method produces optimal solutions in all trials that best represent the theoretical Pareto frontier, with all points falling on the frontier.

## VI. DISCUSSION

This paper investigates and formulates the inspection systems at ports-of-entry. The formulation is general and applicable to different systems as the attributes of a typical container are expressed by a Boolean function. The inspection stations in the physical configuration of the system can be arranged in series (sequential inspection), parallel, series-parallel, parallel-series, *k*-out-of-*n*, or in any network configuration. Boolean functions corresponding to any of these configurations can be developed. The number of attributes and the inspection sequence have significant impact on the system performance. Likewise, the threshold levels of the sensors are critical in the decision process of accepting or classifying a container as suspicious. They influence the probability of making the "wrong" decision by accepting undesired containers or subjecting "good" containers to further unneeded inspections. The POE inspection system problem is formulated as a multi-objective optimization problem that attempts to minimize the total cost as well as the delay time of the containers. The paper presents three different approaches for determining the optimum inspection sequence and the threshold levels at each inspection station that result in the optimization of the system performance measures in terms of cost and time. They are: grid search, constrained nonlinear optimization function, and genetic algorithm. All result in the same values of the optimization function when the number of inspection stations and threshold levels are small. The first two approaches become impractical when more stations and threshold levels are introduced while the genetic algorithm provides optimum or near optimum solutions for such problems in much smaller computation times. As stated earlier, these approaches provide Pareto frontier optimal solutions where every solution includes the optimum sequence of the inspection stations and the corresponding optimum threshold levels. This will enable the decision maker to choose amongst solutions that meet other constraints such as budget, space or layout of the port.

ACKNOWLEDGMENT

This research is conducted with partial support from Office of Naval Research grant numbers N00014-05-1-0237, N00014-07-1-029 and National Science Foundation grant number NSFSES 05-18543 to Rutgers University.